\documentclass[11pt,a4paper,leqno]{article}
\usepackage{latexsym}
\usepackage{amsfonts}
\usepackage {amsbsy}
\usepackage {amsmath}
\usepackage{amssymb}
\newtheorem{defn}{Definition}[section]
\newtheorem{thm}[defn]{Theorem}
\newtheorem{prop}[defn]{Proposition}
\newtheorem{lem}[defn]{Lemma}
\newtheorem{cor}[defn]{Corollary}
\newtheorem{ex}[defn]{Example}

\newcommand \psh{plurisubharmonic }

\newcommand \MA{Monge-Amp\`ere }
\newcommand \demo{ Proof: }
\newcommand \C{\mathbb C}
\newcommand \D{\mathbb D}
\newcommand  \N{\mathbb N}

\newcommand \R{\mathbb R}
\newcommand \B{\mathbb B}
\newcommand \fin{$\blacktriangleright$\\}
\newcommand \lto{\longmapsto}
\newcommand \lra{\longmapsto}

\newcommand \mc{\mathcal}
\newcommand \mb{\mathbb}
\newcommand \mrm{\mathrm}
\newcommand \Om{\Omega}
\newcommand \om{\omega}
\newcommand \sm{\setminus}
\newcommand \es{\emptyset}
\newcommand \sub{\subset}
\newcommand \Sub{\Subset}

\newcommand \de{\delta}

\newcommand \la{\lambda}

\newcommand \ep{\varepsilon}
\newcommand \f{\varphi}

\newcommand \bd{\partial}
\newcommand \we{\wedge}
\numberwithin{equation}{section}
 \title{Partial pluricomplex energy and integrability exponents
 of plurisubharmonic functions }
 \author{ P. {\AA}hag, U. Cegrell, S. Ko\l odziej, H.H. Ph\d{a}m and  A. Zeriahi}
\date{}
\begin{document}
\maketitle

\begin{abstract} We first prove a quantitative estimate of the volume of the sublevel sets of a \psh function in a hyperconvex domain with boundary values $0$ (in a quite general sense) in terms of its Monge-Amp\`ere mass in the domain. Then we deduce a sharp sufficient condition on the Monge-Amp\`ere mass of such a plurisubharmonic function $\f$ for exp$(-2 \f)$ to be globally integrable as well as locally integrable. 
\end{abstract}

\section {Introduction}

It is well known  that estimates on volumes and capacities of sub-level sets of plurisubharmonic  functions from various classes as well as integrability theorems  for such classes play an important role in many areas of complex analysis  (see \cite{Ki2}, \cite{Ko2}, \cite {Ze1}, \cite{Ze2} and references therein). 

A classical result in this direction is Skoda's theorem  \cite{Sk} which asserts that if  $\f$ is a \psh function defined near some point $a \in \C^n$, then $\exp (- 2 \f) $ is locally integrable in a neighbourhood of $a$ if its Lelong number satisfies $\nu_a  (\f) < 1$.  

For $n = 1$ the condition $\nu_a (\f) < 1$ turns out to be equivalent to the local integrability of  $\exp (- 2 \f) $
in a neighbourhood of $a$.
It is possible in this case to derive a global integrability result using classical potential theory. Namely, if $\f$ is a subharmonic function defined
on the unit disc $\D\subset\C$ with smallest harmonic majorant identically zero and
$2 \pi \mu :=\int_{\mathbb D} \Delta \f  <+\infty$. Then for any $s>0$
\begin{equation} \label{intr_est1}
V_2 (\{ \f\leq -s \}) \leq  4 \pi  \exp [ -2 s / \mu  ]\, ,
\end{equation}
where $V_2$ is the $2-$dimensional Lebesgue measure on $\C$. 

From this inequality it is easy to derive a uniform bound on $\int_{\D} e^{- 2 \f} d V_2$ when $\int_{\mathbb D} \Delta \f  \leq 2 \pi \mu < 1$ (see section 4).

When $n \geq 2,$ the situation is much more delicate (see \cite{BJZ} for a partial result).

In \cite{De3}, Demailly provided a sharp condition for the local
integrability of $e^{-2 \f}$, if $\f$ is \psh , in terms of the 
mass of its \MA measure $(dd^c \f)^n .$ It says that if $\Om \Sub \C ^n $, $\f \in PSH(\Om )$ satisfies  $-A \leq \f \leq 0$ on $\Om
\setminus K ,$ where $K\Subset \Om $ and
\begin{eqnarray}
 \int _{\Om } (dd^c \f)^n \leq
\mu ^n <n^n \label{eq:1} \end{eqnarray} then
$$\int _K e^{-2 \f} dV_{2 n} \leq C(\Om , K, A, \mu ) $$
where $dV_{2 n}$ denotes the $2 n-$dimensional Lebesgue measure and  $dd^c = \frac {\sqrt{- 1}}{\pi} \partial \bar{\partial}$. In an appendix to  \cite{De3}, the last-named author observed that from this estimate one can actually deduce a global estimate on the whole of $\Om$.

This result can be viewed as a non linear version of Skoda's integrability
theorem \cite{Sk} where the  assumption that $u$ is bounded near $\bd\Om
$ gives a much stronger statement. Without any extra hypothesis
the estimate is no longer true as functions depending on one
variable only show.

Actually Demailly proved, using an  approximation theorem  \cite{De2}  and the semicontinuity theorem for complex singularity exponents of \psh functions \cite{DK}  (both rather difficult) that his criterion is equivalent to a
local algebra inequality due to Corti  \cite{Co} for $n=2, $ and  de
Fernex, Ein, Musta\c{t}\v{a} \cite{DEM2} in the general case. The inequality
\begin{eqnarray}
lc ({\mc I}) \geq n e({\mc I})^{-1/n} ,\label{eq:2}
\end{eqnarray}
relates $e(\mc I)$ - the Hilbert-Samuel multiplicity of the ideal of
germs of holomorphic functions with isolated singularity at the
origin in $\mb C ^n  $ and $lc(\mc I)$ -  the log canonical threshold
of $\mc I$.

\medskip

Introducing his result Demailly called for an "analytic proof" of
it for the following reasons:

$\bullet $ the criterion involves only \psh\  functions,

$\bullet$ in his proof the constant $C(\Om , K, A, \mu )$ is not
explicitly given in terms of $\Om , K, A, \mu ,$

$\bullet$ an alternative proof combined with Demailly's argument
would provide an analytic way of proving (\ref{eq:2}).

We refer to  \cite{Co},  \cite{De3}, \cite{DEM1} and \cite{DEM2} for the discussion of the
interesting consequences  (\ref{eq:2}) has in the study of
birational rigidity of varieties.

Our first aim is to
generalize~(\ref{intr_est1}) into the several complex variables setting. 
In particular this gives a positive answer to the question of Demailly.
 
\medskip
In order to state our main results let us introduce some notations.
Throughout this paper,  $\Om \Sub \C^n$ denotes a bounded hyperconvex domain (see Section 2 for the definition).  The normalized operator $d^c=\frac{\sqrt{- 1}}{2\pi}(\bar\partial -\partial)$ is used, 
 so that the complex Monge-Amp\`{e}re measure given by $\log |z|$ is exactly the Dirac measure at the origin, i.e. $(dd^c\log|z|)^n=\delta_0$.  

We shall consider the class $\mc E (\Om)$ introduced in~\cite{Ce2}. It is
 essentially the largest set of non-positive plurisubharmonic functions defined on the hyperconvex domain
$\Omega$ for which the complex Monge-Amp\`{e}re operator is
well-defined  (Theorem~4.5 in~\cite{Ce2}). Let $\mc F (\Omega)\subset \mc E (\Omega)$ contain those
functions with smallest maximal plurisubharmonic majorant identically zero and also with finite total Monge-Amp\`{e}re
mass. Note that if $n=1$, then the condition of belonging to $\mc F (\Omega)$ coincides with the  above
conditions on $\f$ that its smallest harmonic majorant is identically zero.

Our main result is the following generalization of ~(\ref{intr_est1}) (see Theorem 4.1 below).
\bigskip

\noindent {\bf Theorem~A} {\it 
There exists a uniform constant $c_n>0$, depending only on $n$  such that for any $\f\in \mc F (\Om)$,
and any $s>0$,  we have that}
\begin{equation}\label{eq:E2}
V_{2 n} (\{ \f \leq -s \})\leq c_n   \delta _\Om^{2 n}  \left(1+s\mu^{-1}\right)^{n-1} \exp \left(-2ns\mu^{-1}\right)\, ,
\end{equation}
{\it where $V_{2 n}$ is the $2 n$-dimensional Lebesgue measure on $\C^n$, $\mu\geq 0$ is defined through
$\mu^n=\int_\Omega(dd^c\f)^n$ and $\delta _\Om$ is the diameter of $\Om$.}

\medskip

As a consequence of our theorem we get a precise and global quantitative version of the integrability theorem of Demailly (see Section 5).

\bigskip

\noindent {\bf Theorem~B} {\it There exist a uniform constant $a_n>0$, depending only on $n$,  such that for any positive number $0 \leq \mu < n$ and any  $\f \in  \mc F (\Omega)$ such that $\int_{\Om} (dd^c \f)^n \leq \mu^n, $ we have that}
\begin{equation}\label{intr_est2}
\int_{\Om} e^{- 2 \f} d V_{2 n}  \leq  \left(\pi^n +  a_n \frac{\mu}{(n - \mu)^n}\right)   \delta _\Om^{2 n},
\end{equation}
{\it where $V_{2 n}$ is the $2 n$-dimensional Lebesgue measure on $\C^n$ and $\delta _\Om$ is the diameter of $\Om$.}

\medskip

The proof of Theorem~A goes by induction on the dimension $n \geq 1$ starting from (\ref{intr_est1}). The first step is to reduce to the case where $\Om = \D^n$  is the unit polydisc using a subextension theorem (\cite{CZ}). Then the  key ingredient in the induction process relies on the special properties of the energy of the slices of  a function $\f \in \mc F (\Om \times D)$ defined on a product domain (see Section 3). To be a bit more precise, we prove the following result   (see Theorem~3.1, Section 3).

\bigskip

\noindent {\bf Theorem~C} {\it 
Let $\Om \Sub \C^n$ and $D \Sub \C$ be hyperconvex domains and $\f \in \mc F
(\Omega \times D)$. Then for almost all $\zeta \in D$, the energy of the slice function  $\f (\cdot,\zeta)$  is well-defined  by} 
\[
\f_{n + 1} (\zeta) :=\int_{\Omega_z}
\f (z,\zeta)\left(dd^c_{z}\f (z,\zeta)\right)^{n}\, .
\] 
 {\it and is equal to 
$$\f_{n + 1} (\zeta)=\int_{\Omega_z \times D_{\eta}}
g(\zeta,\eta)(dd^c \f(z ,\eta))^{n + 1},$$
  where $g = g_{D}$ is the Green
function of $D$ with logarithmic pole. Moreover, $\f_{n + 1} \in \mc F (D)$ and its Laplace mass is given by
$$\int_D dd^c \f_{n + 1} = \int_{\Omega \times D} (dd^c \f)^{n + 1}.$$}
\medskip

In Section ~5 we prove a strong global  version of the integrabilty theorem of Demailly as well as its local version.

In Section 6 we give different applications of our results. We prove a useful inequality between the volume and the \MA capacity of a Borel set (Proposition 6.1) and an integral estimate of \MA capacities of slices of a Borel set in a product domain (Proposition 6.2). Finally in section 6.3 we deduce a general local transcendental  inequality on complex singularity exponents of \psh functions (Proposition 6.3) which implies,  following an argument of Demailly in \cite{De3}, the local algebra inequality (\ref{eq:2}).

\section {Preliminaries} 
 Let us  recall some definitions. Let $\D$ denote the unit disk in $\C$, $\D^n$ the unit poldisc in $\C^n$ and let $V_{2 n}$ denote the Lebesgue measure in $\C^n $. Consider also 
  the usual differential operators $d$ and  $d^c =\frac{\sqrt{-1}}{2\pi }(\bar{\partial} -\partial )$  acting on \psh functions on domains in $\C^n$ so that $dd^c = (\sqrt{- 1} \slash \pi ) \partial \bar \partial$.

For an open set $\Om \sub \C^n$, we denote by   $PSH (\Om) \sub L^1_{loc} (\Om)$ the set of plurisubharmonic functions in $\Om$.  

An open set $\Om \Sub \C^n$ is said to be hyperconvex if it admits a negative plurisubharmonic  exhaution function i.e. there exists a \psh function $\rho : \Om \lra [- 1, 0[$ such that for any $c < 0$, $\Om_c := z \in \Om ; \rho (z) < c\} \Sub \Om.$  

It is well known that a domain $D \Sub \C$ is hyperconvex if and only if it is regular with respect to the  Dirichlet problem for the Laplace operator (\cite{Ra}).  Therefore any product of regular planar domains (e.g. a polydisc) is a hyperconvex domain. More generally any bounded pseudoconvex domain with Lipschitz boundary is hyperconvex (see \cite{De1} and references therein).

Now we recall some notations from \cite{Ce1}, \cite{Ce2}. We write $\mc E_{0} (\Om )$ for the set of \psh test functions i.e.  functions $ \f \in PSH (\Om) \cap L^{\infty} (\Om)$ which tend to zero at the boundary
and satisfy $\int_{\Om } (dd^c \f )^n < + \infty.$

 Denote by $\mc F (\Om )$ the set of all
$\f \in PSH (\Om )$ such that there exists a sequence $(\f _{j})$
of \psh functions in $\mc E_{0} (\Om )$ such that $\f _{j}
\searrow \f$ and $\sup_{j} \int_{\Om } (dd^c \f _{j})^n < +
\infty.$

 The class $\mc E (\Om )$ will be the set of all
$\f \in PSH (\Om )$ such that for any open subset $\om \Subset \Om$ there is a function $\psi \in \mc F (\Om )$ such that  $\psi = \f$ on $\om$.

The complex Monge-Amp\`ere operator is well defined and continuous under decreasing limits in the class $\mc E (\Om )$. Moreover
in the class $\mc F (\Om )$, we have the following strong convergence theorem, namely if $(\f_j)$ is a decreasing sequence of functions in $\mc F (\Om )$ which converges to $\f \in \mc F (\Om )$, then for any $h \in PSH (\Om)$ such that $h \leq 0,$ we have (see \cite{Ce2}, \cite{Ce4})
$$\lim_j \int_{\Om} h (dd^c \f_j)^n = \int_{\Om} h (dd^c \f)^n.$$
 
Define $\mc E_1 (\Om)$ to be the class of plurisubharmonic functions 
$\f \in PSH (\Om )$ with finite energy i.e. there exists a sequence $(\f _{j})$
of \psh functions in $\mc E_{0} (\Om )$ such that $\f _{j}
\searrow \f$ and $\sup_{j} \int_{\Om } (-\f_j) (dd^c \f _{j})^n < +
\infty.$ It can be proved that $ \mc E_1 (\Om) \sub \mc E (\Om)$ (see \cite{Ce2}).

We will need the following lemma.
\begin{lem} Let $ v \in PSH (\Om) \cap L^{\infty} (\Om)$ be such that $\lim_{z \to \zeta} v (\zeta) = 0$ for any $\zeta \in \partial \Om$. Assume 
that 
$\int_{\Om} ( - v) (dd^c v)^n < + \infty.$ Then $ v \in \mc E_1 (\Om)$.
\end{lem}
\demo Let $(\Om_j)$ be an exhaustion of $\Om$ by bounded domains. It follows from 
\cite{Ko1} that for $j \in \N,$ there exists $v_j \in \mc E_0 (\Om)$ such that
$$
(dd^c v_j)^n = {\bf 1}_{\Om_j}(dd^c v)^n $$
in $\Om$.
By the comparison principle $(v_j)_j$ is a decreasing sequence from $\mc E_0 (\Om)$ converging to $v$. Integration by parts gives that
$\int_{\Om} (- v_j) (dd^c v_j)^n \leq  \int_{\Om} (- v) (dd^c v)^n < + \infty,$ so $v \in \mc E_1 (\Om)$ by definition.
\fin

Now we introduce the notion of capacity due to Bedford and Taylor (\cite{BT}).
For a given Borel subset $E \sub \Om$ we define the Monge-Amp\`ere capacity of the condenser $(E,\Om)$ by
$$
Cap (E,\Om) = Cap_{\Om} (E) := \sup \{ \int_E (dd^c v)^n ; v \in PSH (\Om), -1 \leq v \leq 0\}.
$$
Then by \cite{BT} if $E \Sub \Om$ we have the formula
$$
Cap (E,\Om) = \int_{\Om} (dd^c h^*_{E,\Om})^n,
$$
where $h_{E,\Om}$ is the extremal function of $(E,\Om)$ defined by
$$
h_{E,\Om} := \sup \{v \in PSH (\Om) ; v  \leq 0, v|E \leq - 1\}.
$$
We will also need the following estimates on the capacity of the sub-level sets of functions in  $\mc E_1 (\Om)$ (see \cite{CKZ}).

\begin{lem} Let  $v \in \mc E_1 (\Om)$. Then for any $s > 0$, we have that
$$
s^{n +1} Cap_{\Om} (\{v \leq- s\}) \leq \int_{\Om} (- v) (dd^c v)^n.
$$
\end{lem}
\demo
By homogeneity, it is enough to prove the estimate for $s = 1$. 
Then take an arbitrary compact subset $K \sub \{v \leq- 1\}$. If $h_K$ is the extremal function of $(K,\Om)$ the function
 $h :=  h_K^* \in \mc E_0 (\Om)$ and satisfies $ v \leq h$. Thus using 
repeatedly  integration by parts we obtain  that
\begin{eqnarray}
 Cap (K,\Om) &=&\int _{\Om } (-h)(dd^c h)^n \leq \int _{\Om} (-v) (dd^c h)^n \nonumber \\ 
& \leq & \int _{\Om } (-h )dd^c
v \we (dd^c h )^{n-1} \leq \int _{\Om} (-v )dd^c v \we (dd^c h
)^{n-1} \\ & \leq & ... 
  \leq \int _{\Om} (-v
)(dd^c v )^n .\nonumber
\end{eqnarray} 

\section{Partial energies}

For a hyperconvex domain $D \sub \C$, we denote by $g = g_{D}$ the Green function of $D$ with logarithmic pole.

We prove the following result (Theorem C in the introduction).
\begin{thm}
Let $\Omega \subset\C^{n}$ and  $D \subset \C$ be
two bounded hyperconvex domains. If $\f \in \mc F
(\Omega \times D)$, the slice function $\Omega  \ni
z \to \f (z,\zeta)\in\mc E_1 (\Omega)$
 for all $\zeta \in D$ with
\[
\int_{\Omega_z \times D_{\eta}}
g(\zeta,\eta)(dd^c(\f (z,\eta))^{n + 1} > -\infty. 
\]
 Furthermore, if
we define
\[
\f_{n + 1} (\zeta) :=\int_{\Omega_z}
\f (z,\zeta)\left(dd^c_{z} \f (z,\zeta)\right)^{n}\,\]
if the integral is well-defined and $\f_{n + 1} (\zeta) = -\infty$ otherwise,
then for any $\zeta \in D,$ 
$$\f_{n + 1} (\zeta) = \int_{\Omega_z \times D_{\eta}}
g(\zeta,\eta)(dd^c \f (z,\eta))^{n + 1}.$$
 In particular we have that
$\f_{n + 1} \in\mc F (D)$ and it satisfies
$$\int_D
dd^c \f_{n + 1} = \int_{\Omega \times D}(dd^c \f)^{n + 1}.$$
\end{thm}
\demo Assume first that $ \f \in \mc E_0 (\Omega \times D) \cap C^{\infty}(\Omega \times D)$, and let
$K \Subset \Omega$, $L \Subset D$. Then $0\geq
\f (z,\zeta)\left(dd^c_{z}\f  (z,\zeta)\right)^{n}\in
C^{\infty}(\Omega \times D)$ and thus
\[
\f^K (\zeta) :=\int_{K}
\f (z,\zeta) \left(dd^c_{z}\f (z,\zeta) \right)^{n}\in
C^{\infty}(D)\, .
\]
For  $h \in \mc E_0 (D)\cap C(D)$ we have 
\begin{eqnarray*}
\int_{L} \f^K(\zeta) dd^c h(\zeta) & =&\int_{L} \int_{K}
\f (z,\zeta) \left(dd^c_{z}\f (z,\zeta)\right)^{n}\wedge
dd^ch(\zeta)\\
& =& \int_{L} \int_{K}
\f (z,\zeta)\left(dd^c \f (z,\zeta)\right)^{n}\wedge
dd^ch(\zeta).
\end{eqnarray*}
Then it follows that  
\begin{eqnarray*}
&\int_{L} \int_{K}
 \f (z,\zeta) \left(dd^c \f (z,\zeta)\right)^{n}\wedge
dd^ch(\zeta) \geq \\
& \int_{D} \int_{\Om}
\f (z,\zeta)\left(dd^c \f (z,\zeta)\right)^{n}\wedge
dd^ch(\zeta).
\end{eqnarray*}
By a generalized Jensen-Lelong-Demailly type formula in $\mc F (\Om \times D)$ (\cite{CK}, Remark 1), we have
\begin{eqnarray*}
& \int_{\Om \times D}
\f (z,\zeta) \left(dd^c \f (z,\zeta)\right)^{n}\wedge
dd^ch(\zeta) = \\
& \int_{\Om \times D}
h (\zeta) \left(dd^c \f (z,\zeta)\right)^{n + 1} > - \infty,
\end{eqnarray*} 
since $h$ as a function of $(z,\zeta) \in \Om \times D$ which only depends on $\zeta \in D$ and vanishes on the distinguished boundary of $\Om \times D$.

Then by letting $L$ increase to $D$, it follows that $u^K$ is a decreasing family of continuous functions on $\Om$ which are uniformly integrable
 on $\Om$ as 
 $K$ increases to $\Om$. This implies that $\f_{n + 1}$ is upper semi-continous and integrable on $\Om$ and satisfies
\begin{equation}
\int_{D} \f_{n + 1} dd^c h = \int_{D_{\zeta}} \int_{\Om_z}
h (\zeta)\left(dd^c \f (z,\zeta)\right)^{n},
\label{eq:F}
\end{equation}
for any test function $h \in \mc E_0 (D) \cap C(D)$.
Since $C_{0}^{\infty} (D)\subset \mc E_0(D)\cap
C(\bar{D}) -  \mc E_0(D)\cap C(\bar D)$
(see \cite{Ce2}), we get from (\ref{eq:F}) that $dd^c \f_{n + 1} \geq 0$ in the weak sense on $\Om$, which proves that $\f_{n + 1}$ is subharmonic on $D$.

Now fix $\zeta \in D$ and apply (\ref{eq:F}) to the function $h = \sup \{g (\zeta,\cdot),-j\}$. Then by classical potential theory in $\C$, we deduce that
\begin{equation}
\label{eq}
 \f_{n + 1} (\zeta)=\int_{\Omega_z \times D_{\eta}}
g(\zeta,\eta) (dd^c \f (z,\eta))^{n+1},
\end{equation}
since $dd^c g(\zeta,\cdot)$ is the Dirac mass at the point $\zeta$.
This also proves that $\f_{n + 1} \in\mc F (D).$

If $\f_{n + 1} (\zeta)>-\infty$, then
$v := \f (\cdot,\zeta)$ has boundary values $0$ and
$$
\int_{\Om} (-v) (dd^c v)^{n + 1}  = -  \f_{n + 1} (\zeta) < + \infty,
$$
 which implies by Lemma 2.1 that $v =  \f (\cdot,\zeta) \in\mc E_1(\Omega)$.

 For the general case, assume that $\f \in\mc F(\Omega \times D)$. By \cite{Ce3} we can choose a sequence
 $\f^j\in\mc E_0\cap C^{\infty}(\Omega \times D)$ such that $\f^j\searrow \f$, $j\to +\infty$. It follows
 from~(\ref{eq}) that $\left(\f^j_{n +1}\right)_j$ is a decreasing sequence of functions in $\mc F (D)$ such that
 \[
 \lim_{j\to +\infty}\f^j_{n+1}(\zeta) =\int_{\Omega \times D}g(\zeta,\eta) (dd^c{\f})^{n +1}\, .
 \]
 It follows now from the previous case that $\left(\f^j(\cdot,\zeta)\right)_j$ is a decreasing sequence 
of functions in  $\mc E_1(\Omega)$ with uniformly 
bounded energies  which converges to $\f (\cdot,\zeta)$ if $\f_{n + 1} (\zeta)>-\infty$. 
Then from   Theorem 3.8 in  \cite{Ce1}, we deduce that
 $\f (\cdot,\zeta)\in\mc E_1(\Omega)$ if $\f_{n + 1} (\zeta)>-\infty$.  
\fin

\begin{ex} The function
$$ \f (z,\zeta) := \sum_{k = 1}^{+ \infty} \max \{\log \vert z\vert, k^{- 4} \log \vert \zeta \vert\}, \ (z,\zeta) \in \D \times \D$$
is an example of a function $\f  \in \mc F (\D^2)$ with all  slices  $\f (\cdot,\zeta) \in \mc E_1 (\D) \sm \mc F (\D)$ if $\zeta \neq 0$ (see Example 5.7 \cite{CW}).
\end{ex}

The last result can be generalized as follows.

\begin{thm}
Let $\Omega \subset\C^{n}$ and  $D \subset\C$ be two bounded hyperconvex domains.
If $\f_j \in \mc F (\Omega \times D), 0\leq j \leq n$,  the slice function $\Omega \ni z \to \f_0(z,\zeta) \in L^1(\Omega)$
with respect to $dd^c_{z} \f_1(z,\zeta) \wedge\dots\wedge dd^c_{z} \f_n (z,\zeta)$ for all $\zeta \in D$ with
$$
\int_{\Omega_z \times D_{\eta}} g(\zeta,\eta) dd^c \f_0(z,\eta) \wedge dd^c_{z}\f_1(z,\zeta) \wedge \dots \wedge dd^c \f_n(w,\eta)>-\infty\,.
$$
 Furthermore, if we define
\[
u (\zeta) =\int_{\Omega_z} \f_0(z,\zeta) \  dd^c_{z} \f_1(z,\zeta) \wedge\dots\wedge dd^c_{z} \f_n(z,\zeta) \, .
\]
then 
$$
u (\zeta)=\int_{\Omega_z \times D_{\eta}} g(\zeta,\eta) \ dd^c \f_0(z,\eta) \wedge dd^c \f_1(z,\eta) \wedge\dots \wedge dd^c \f_n(z,\eta)\, 
$$
so in particular we have that $u \in \mc F (\Omega)$ and 
$$
 \int_\Omega dd^c u=\int_{\Omega \times D} dd^c \f_0 \wedge dd^c \f_1 \wedge\dots\wedge dd^c \f_n.
$$
\end{thm}
\demo
Let $h\in\mc E_0 \cap C(D)$ be a given test function.
As in the first part of the proof of Theorem 3.1 we get
$$
\int_D u \ dd^c h = \int_{\Omega \times D} \f_0(z,\zeta) \ dd^c \f_1(z,\zeta)\wedge\dots\wedge dd^c \f_n(z,\zeta)\wedge dd^c h (\zeta). 
$$
So if we prove that the right hand side equals
$$
\int_{\Omega \times D} h(\zeta) \  dd^c \f_0(z,\zeta)\wedge\dots\wedge dd^c \f_n(z,\zeta), 
$$
the proof can be completed in the same way as in the second part of the proof of the previous theorem.

Indeed for $0\leq k_j \leq n, 0\leq j \leq n$ the inequality
$$
\int_{\Omega \times D} \f_{k_0}(z,\zeta) \ dd^c \f_{k_1}(z,\zeta)\wedge\dots\wedge dd^c \f_{k_n}(z,\zeta)\wedge dd^ch(\zeta) \geq 
$$
$$
\int_{\Omega \times D} h(\zeta) \  dd^c \f_{k_0}(z,\zeta) \wedge dd^c \f_{k_1}(z,\zeta) \wedge\dots\wedge dd^c \f_{k_n}(z,\zeta), 
$$
can be obtained by approximating $h$ by a decreasing sequence of functions $h_p\in\mathcal E_0(\Omega \times D),$
using partial integration in $\mathcal F$ and observing that 
$dd^c \f_{k_1}(z,\zeta) \wedge\dots\wedge dd^c \f_{k_n}(z,\zeta)\wedge dd^ch_p(\zeta )$ tends weakly to
$dd^ c \f_{k_1}(z,\zeta) \wedge\dots\wedge dd^c \f_{k_n}(z,\zeta)\wedge dd^c h(\zeta)$ 
when $p$ tends to $\infty.$

Again, since $\mathcal F$ is a convex cone, $\psi := \sum_{j=0}^n \f_j \in \mc F (\Omega \times D)$ and then  it follows from \cite{CK} that
$$
\int_{\Omega \times D} \psi (z,\zeta) \ \left(dd^c \psi (z,\zeta)\right)^{n}\wedge dd^ch (\zeta) =
$$
$$
\int_{\Omega \times D} h(\zeta) \  \left(dd^c \psi (z,\zeta)\right)^{n+1}.
$$

Using the separate linearity of the wedge product, we can expand both sides to obtain sums of terms of the form 
$$
\int_{\Omega \times D} \f_{k_0}(z,\zeta) \ dd^c \f_{k_1}(z,\zeta) \wedge \dots \wedge dd^c \f_{k_n}(z,\zeta) \wedge dd^ch (\zeta)
$$
on the left hand side, while on the right hand side we get terms
$$
\int_{\Omega \times D} h (\zeta) dd^c \f_{k_0} (z,\zeta) \wedge dd^c \f_{k_1}(z,\zeta)\wedge \dots \wedge dd^c \f_{k_n} (z,\zeta), 
$$
Since they have the same sum, they have all to be equal which completes the proof of the theorem.
\fin

\section{Volume estimates of sub-level sets}

Here we prove the following result.

\begin{thm}  There exists a constant $c_n > 0$ such that for any $\mu >  0$, any $ \f \in
\mc F (\D^n) $ with $\int_{\D^n} (dd^c \f)^n \leq \mu^n$ and any  $s > 0$, we have that
\begin{eqnarray}
V_{2 n}  \left(\{\f \leq - s\}\right) \leq  c_n (1 + s \slash \mu)^{n - 1}  \exp \left( - 2 n s \slash \mu \right).
\label{eq:VMA} 
\end{eqnarray}
\end{thm}
\noindent{\demo} We prove the theorem  using induction over the dimension $n$. 

 For $n=1$, the estimate was proved in (\cite{BJZ}). Let us recall it here for the convenience of the reader.  We use the classical  P\'{o}lya's inequality  which we recall.  Let $K \sub \C$ be a compact subset in the complex plane with area $A (K)$ and logarithmic capacity  $c(K)$. In ~\cite{P}, P\'{o}lya proved what we today could write as 
\begin{equation}
A (K) \ \leq \ \pi  c(K)^2\, ,
\label{eq:P}
\end{equation}
 (for an elegant proof see e.g. Theorem~5.3.5 in~\cite{Ra}). 
Now assume that $K$ is not polar. Then from the Riesz representation formula for the Green function $V_K$  of $K$ with pole at infinity, we obtain that
$$ c (K) \leq 2 \exp \left(- \sup_{\vert z\vert = 1} V_K (z)\right).$$
Now denote by $M_K := \sup_{\vert z\vert = 1} V_K (z).$ Then $M_K^{- 1} V_K \leq h_{K^,\D}$ on $\D$.
Since by our normalization $ \int_{\D}
 dd^c V_K =  \int_{\C}
 dd^c V_K = 1$, it follows from the comparison principle for the Laplace operator that
$$ M_K^{- 1} = M_K^{- 1} \int_{\D}
 dd^c V_K \leq \int_{\D} dd^c h_K^* = Cap (K,\D).$$
This inequality is due to Alexander and Taylor (see \cite{AT}).
Then putting all together we obtain the inequality  
$$V_2 (K) \leq 4 \pi \exp \left(- 2 \slash Cap (K,\D)\right).$$
It is clear that this inequality is still true for Borel subsets $K \sub \D$. 
Now if $\f \in \mc F (\D),$ we know that $Cap (\{\f \leq - s\},\D) \leq \int_{\D} dd^c \f  \slash s$. Therefore we have that
\begin{equation}
V_2 \left(\{\f \leq - s\}\right) \  \leq 4 \pi  \  \exp \left( -2 s \slash \mu\right),
\label{eq:P}
\end{equation}
where $\mu := \int_{\D} dd^c \f$  (see \cite{BJZ}).

Now assume that the estimate (\ref{eq:VMA}) is true in dimension $n$ and let us prove it in dimension $n + 1$.

 Fix   $\f \in \mc F (\D^{n + 1})$ such that 
$$
\int _{\D ^{n+1} } (dd^c \f)^{n+1} \leq \mu^{n + 1}.
$$
By homogeneity, it is enough to prove the estimate for $s = 1.$ Then we want to estimate the volume $V_{2 n + 2} \left(\{\f \leq - 1\}\right)$ by applying Fubini's Theorem. So fix $\zeta \in \D$ and  estimate the volume $V_{2 n} (\{ z \in
 \D ^n : \f (z , \zeta) \leq -1 \})$. Indeed, define   $E_{\zeta} := \{ z \in
 \D ^n : \f (z , \zeta) \leq -1 \}$, consider its relative extremal function $h_{\zeta} := h_{E_{\zeta}}^*$ and observe that
$V_{2 n} (E_{\zeta}) =  V_{2 n} (\{h_{\zeta} \leq - 1\})$, since the two sets coincide up to a pluripolar set.
We want to apply the induction hypothesis to the function $h_{\zeta}$. 

Fix $\zeta \in \D$ such that $\f_{n + 1} (\zeta) > - \infty$ and 
observe  that  $ h_{\zeta}  \geq  v := \f (\cdot , \zeta)$. By Theorem 3.1, the function $v = \f (\cdot , \zeta) \in \mc E_1(\D^n)$ and then
 $  h_{\zeta} \in \mc E_1(\D^n)$.
On the other hand, by \cite{BT}, we know that
$$\int_{\D^n} (dd^c h_{\zeta})^n =  Cap (E_{\zeta},D^n).$$
Then since $v = \f (\cdot , \zeta) \in \mc E_1(\D^n)$, it follows from Lemma 2.2 that
$$
Cap (E_{\zeta},D^n) = Cap (\{v \leq - 1\})  \leq \int_{\D^n} (- v) (dd^c v)^n = - \f_{n+1} (\zeta) < + \infty.
$$

This implies that for any $\zeta \in \D$ such that $\f_{n + 1} (\zeta) > - \infty$,
$$  \int_{\D^n} (dd^c  h_{\zeta})^n  \leq  - \f_{n +1} (\zeta) < + \infty$$
and then $h_{\zeta} \in \mc F (\D^n)$.

Now applying the induction hypothesis to the function $h_{\zeta} \in \mc F (\D^n)$, we deduce that
for almost all $\zeta \in \D$,
$$
V_{2 n} \left(\{ \f(\cdot,\zeta) \leq - 1\}\right) \leq c_n \left( 1 + (- \f_{n +1} (\zeta))^{- 1 \slash n}\right)^{n - 1} 
\exp \left(- 2 n (- \f_{n +1} (\zeta))^{- 1 \slash n}\right).
$$
Then integrating in $\zeta \in \D,$ we get
\begin{equation}
V_{2 n + 2} \left(\{ \f  \leq - 1\}\right) 
\leq  c_n \int_{\D} \chi \left(- \f_{n + 1} (\zeta)\right) d V_2 (\zeta),
\label{eq:I}
\end{equation}
where  
$$\chi (t) := \left(1 + t^{- 1 \slash n}\right)^{n - 1} \exp \left(- 2 n  t^{- 1 \slash n}\right), t \geq 0.$$

It is easy to check that  the function $\chi$ is increasing with $\chi (0) = 0$ and $\chi (+ \infty) = + \infty$. Therefore from (\ref{eq:I}), it follows that
\begin{equation}
V_{2 n + 2} \left(\{ \f \leq - 1\}\right)  \leq  c_n \int_0^{+ \infty} \chi'(t) V_2 \left(\{\f_{n + 1} \leq - t\}\right) d t.
\label{eq:II} 
\end{equation}
Since $\int_{\D} dd^c \f_{n + 1} = \int_{\D^{n + 1}} (dd^c \f)^{n + 1} \leq \mu^{n + 1}$ by Theorem 3.1, it follows from 
(\ref{eq:II}) that
$$
V_{2 n + 2} \left(\{ \f \leq - 1\}\right)  \leq  c_n \int_0^{+ \infty} \chi'(t) \exp \left(- 2  t \mu^{- n -1}\right) d t.
$$
Now using the change of variable $x = t^{- 1 \slash n}$ and observing that 
$$ 
\chi' (t) \ d t = -  (2 n x + n + 1) (1 + x)^{n - 2}  e^{- 2 n x} \ d x,
$$
we get the following estimate

\begin{equation}
V_{2 n + 2} \left(\{ \f \leq - 1\}\right)  \leq  8 n \pi  c_n \int_{0}^{+ \infty} (x + 1)^{n - 1} \exp \left( - 2 (n x +  x^{- n}\mu^{- n - 1})\right) d x.
\label{eq:IV}
\end{equation}

Now observe that the function $\R^+ \ni x \lto 2 (n x + x^{-n} \mu^{- n - 1})$ reaches its minimum at the point $ x = 1 \slash \mu$ and 
this minimum is precisely equal to $ 2 (n + 1) \mu^{- 1}$.
Then splitting the integral in (\ref{eq:IV}) into two parts, integrating first from $0$ to $3 \slash \mu$ and then from $3 \slash \mu$ to $+ \infty$, we easily get 
\begin{eqnarray}  
V_{2 n + 2} \left(\{ \f \leq - 1\}\right) 
 &\leq&  8 n \pi c_n  (3 \slash \mu) (1 + 3\slash \mu)^{n -1}  \exp \left(-  2 (n + 1) \mu^{- 1})\right) \nonumber \\
&+& 8 n  \pi c_n \int_{3 \slash \mu}^{+ \infty} (1 + x)^{n - 1} \exp ( - 2 n x ) d x.
\end{eqnarray}
It is easy to see that the last terms is much better that the first one and can be easily estimated from above by  
 $8 \pi  c_n  \exp \left( - 2 (n+1) \slash \mu\right)$ 
 so that we finally get
$$
V_{2 n + 2} \left(\{ \f \leq - 1\}\right)  \leq c_{n + 1} \left(1 + 1 \slash \mu \right)^n \exp \left(-  2 (n + 1) \mu^{- 1})\right),
$$
where $c_{n + 1}:=  8 \pi (n 3^{n} + 1)  c_n$ which implies that 
\begin{equation}
c_{n}:=  2^{3 n - 1}  \pi^n \Pi_{0 \leq k \leq n - 1} (k 3^{k} + 1), \ \ n \geq 1.
\label{eq:c_n}
\end{equation} 
\fin

The volume estimate actually holds in the following general setting which will prove Theorem A.

\begin{cor} Let $\Om \Sub \C^n$ be a bounded hyperconvex domain. Then for any $\f \in \mc F (\Om)$  and any $s > 0$, we have
\begin{equation}
V_{2 n } \left(\{ \f  \leq - s\}\right)  \leq c_n \de_{\Om}^{2 n} \left(1 +  s \mu^{- 1} \right)^{n - 1} \exp \left(-  2 n s \mu^{- 1})\right),
\label{eq:V-Est}
\end{equation}
where  $ \mu^n := \int_{\Om} (dd^c \f)^n$, $\de_{\Om}$ is the diameter of $\Om$  and $c_n $ is the constant defined by (\ref{eq:c_n}).
\end{cor}
\demo Observe that  the inequality is invariant under holomorphic linear change of variables so that we can always assume that $\Om \sub \D^n.$ Then by the subextension theorem (\cite{CZ}), there exists a function $\psi \in \mc F (\D^n)$ such that $\psi \leq \f$ and $\int_{\D^n} (dd^c \psi)^n \leq \int_{\D^n} (dd^c \f)^n = \mu^n.$ Applying the estimate of Theorem 4.1 to $\psi$, we obtain the required estimate.\fin

\section{Integrability theorems in terms of \MA masses}

In this section we prove Theorem B (see Corollary 5.2) stated in the introduction, which will give a pluripotential proof of a theorem due to Demailly \cite{De3}. We also prove a theorem on local integrability.

\subsection{Global integrability}

\begin{thm} Let  $\f \in \mc F (\D^n)$ such that $\int_{\D^n} (dd^c \f)^n \leq \mu^n $ with $\mu < n$. Then 
$$
\int_{\D^n}  e^{- 2 \f} d V_{2 n} \leq \pi^n +   a_n \frac{\mu}{(n - \mu)^n},
$$
 where $a_n > 0$ is a dimensional constant.
\end{thm}

\demo
By Theorem 4.1, we have
\begin{eqnarray*}
\int_{\D^n} e^{- 2 \f} d V_{2 n} &=& \pi^n + 
2 \int_{0}^{+ \infty} e^{2 s} V_{2 n} (\{\f < - s \}) d s \\
&\leq& \pi^n + 2 c_n  \int_{0}^{+ \infty} (1 + s \slash \mu)^{n - 1} e^{2 s - 2 n s\slash \mu} d s.
\end{eqnarray*}

Now it is easy to see by integration by parts that the integrals 
$I_n := \int_{0}^{+ \infty} (1 + s \slash \mu)^{n - 1} e^{2 s - 2 n s \slash \mu} d s$ satisfy the  inequality
$$I_{n } \leq \frac{(n - 1)!}{2^n} \frac{\mu}{(n - \mu)^n}, $$
for $n \geq 1$,
and the required estimate follows with the constant 
\begin{equation}
a_n := \frac{(n - 1)!}{2^{n- 1}} c_n.
\label{eq:a_n}
\end{equation}
\fin
We have a more general result.
\begin{cor} Let $\Om \Sub \C^n$ be a bounded hyperconvex domain.Then  for any $\f \in \mc F (\Om)$ such that $\int_{\Om} (dd^c \f)^n \leq \mu^n < n^n$ with $\mu < n$,  we have that
$$
\int_{\Om} e^{- 2 \f} d V_{2 n} \leq \left( \pi^n + a_n  \frac{\mu}{(n - \mu)^n}\right) \de_\Om^{2 n},
$$
where  $\de_\Om$ is the diameter of  the domain $\Om$ and $a_n$ is the constant defined by (\ref{eq:a_n}).
\end{cor}
\demo
The proof is the same as before using Corollary 4.2.
\fin

As a corollary we get a strengthened version of Demailly's theorem \cite{De2}.

\begin{cor} Let $\Om \Sub \C^n$ be a bounded pseudoconvex domain and $M > 0$ a fixed constant. Then for any $\f \in PSH (\Om)$ with 
$0 \geq \f \geq -  M$  near the boundary, $\f \in \mc E(\Om).$
Moreover if  $\int_{\Om} (dd^c \f)^n \leq \mu^n$ with $\mu  < n$, we have that
$$ \int_{\Om} e^{- 2 \f} d V_{2 n} \leq  \left(\pi^n +  a_n  \frac{\mu}{(n - \mu)^n}\right) e^{2 M}  \de_\Om^{2 n},
$$
where  $\delta _\Om := diam (\Om)$ is the diameter of $\Om$ and $a_n$ is the constant defined by (\ref{eq:a_n}).
\end{cor}
\demo We can assume the domain $\Om$ to be hyperconvex.
It follows from Theorem 2.1 in \cite{Ce4} that there exists $\psi \in \mc F (\Om)$ with $\int_{\Om}(dd^c \psi)^n \leq 
 \int_{\Om}(dd^c \f)^n$  such that $\f \geq \psi - M$ on $\Om$. Then the result  follows now from Corollary 5.2.
\fin

Now we investigate integrability in the critical case when the total \MA mass has the maximal value $n^n$.

 \begin{thm} Let $\Om \Sub \C^n$ be a bounded hyperconvex  
 domain and $\f \in \mc F (\Om)$ such that $\int_{\Om}  (dd^c \f)^n = n^n$. Then for any real number $\la > n$,  we have  that
 $$
 \int_{\Om} \frac{ e^{- 2 \f}}{(1 - \f)^{\la}} d V_{2 n} \leq     
 \left( 1 + (2 \slash \la)^{\la} e^{\la - 2}\right) V_{2 
 n} (\Om) +  2 c_n  \de_{\Om}^{2 n} \frac{1}{\la - n},
 $$
where $c_n$ is the constant defined by (\ref{eq:c_n}).
 \end{thm}
 \demo Indeed, set $\chi (t) := (1 + t)^{- \la} e^{2 t},$ for $t \geq 0$. Since 
 $$
 \chi' (t) = \left(- \la (1 + t)^{- \la - 1} + 2 (1 + t)^{-  
 \la}\right) e^{ 2 t} = (2 - \la  + 2 t) (1 + t)^{- \la - 1} e^{2  t},
 $$
 it follows that the function $\chi$ is increasing for $t \geq t_0 
 := (\la - 2)\slash 2$ and decreasing on $[0,  t_0]$. 

 Therefore we have that
 \begin{eqnarray*}
 \int_{\Om} \frac{ e^{- 2 \f}}{(1 - \f)^{\la}} d V_{2 n} 
 & = & \int_{- \f <  t_0} \frac{ e^{- 2 \f}}{(1 - \f)^{\la}} d 
 V_{2 n} \\
 & + & \int_{- \f \geq t_0} \frac{ e^{- 2 \f}}{(1 - \f)^{\la}} d 
 V_{2 n} \\
 & \leq &
 V_{2 n} (\Om) + \int_{\f \leq - t_0} \frac{ e^{- 2 \f}}{(1 -  
 \f)^{\la}} d V_{2 n} \\
 & \leq &
  V_{2 n} (\Om) + \chi (t_0) V_{2 n} (\Om) + 
  \int_{t_0}^{+ \infty} \chi'(t) V_{2 n} (\{\f \leq - t\}) d t.
  \end{eqnarray*}

 By Corollary 4.2, we have that
 $$
 \int_{\Om} \frac{ e^{- 2 \f}}{(1 - \f)^{\la}} d V_{2 n} 
 \leq \left(1 + \chi (t_0)\right)  V_{2 n} (\Om) + 2 c_n  
 \de_{\Om}^{2 n} \int_0^{+ \infty} (1 + t)^{n - \la - 1} d t.$$
\fin
\subsection{Local integrability of $ exp(-2 \f)$}

\begin{thm}
Suppose $ \f \in\mathcal E(\Omega)$ and $a\in \Omega$.  If $\int\limits_{\{a\}}(dd^c \f)^n<n^n$, then $\exp (-2 \f)$ is locally integrable near $ a$, \\
\end{thm}
\demo We can assume that $\f \in \mc F (\Om).$ Set for $j \geq 1$, 
$$
\psi_j := \sup \{ u \in PSH (\Om) ; u \leq 0, u \leq \f \ \ \text{on} \ \ B_j\},
$$
where $B_j := \B (a,1 \slash j)$ is the ball of center $a$ and radius $1 \slash j$.

Then $\psi_j \in \mc F (\Om), \psi_j \geq \f$  and $\psi_j = \f$ on $B_j$. 
Moreover, $supp (dd^c\psi_j)^n \subset\subset B_{j-1}.$
Denote by $G(z,a)$ the pluricomplex Green function for $\Omega$ with logarithmic pole at $a$ and choose $\delta >0$ so small that 
$$ \int_{\Om} (- \max\{\delta G(z,a),-1\}) (dd^c\varphi)^n < n^n.$$
 Using integration by parts in $\mathcal F (\Om)$ we see that
$$\int_{\Om} \left( - \max\{ \delta G(z,a),-1\}\right)  (dd^c \psi_j)^n \leq \int_{\Om} \left( - \max\{ \delta G(z,a),-1\}\right)  (dd^c\varphi)^n  < n^n. $$
If we choose $ k$ so large that $B_{k- 1} \subset \subset \{\delta G(z,a) < -1\},$ it follows that $\int_{\Omega}(dd^c\psi_k)^n = \int_{B_{k - 1}} (dd^c\psi_k)^n < n^n.$

Now since $\psi_k = \f$ on $B_k,$ it follows from Corollary 5.2 that
$$
\int_{B_k} e^{- 2 \f } d V_n = \int_{B_k} e^{- 2 \psi_k } d V_n \leq \int_{\Om} e^{- 2 \psi_k } d V_n < + \infty.
$$
\fin

\noindent {\bf Remark:} Note that the theorem is optimal as the functions \\
$(n-\ep) \log \vert z - a\vert$ ($\ep > 0$) show.

\section{Applications}

\subsection{An inequality between volume and capacity}
Our first application of Theorem 4.1 is  a useful inequality between volume and Monge-Amp\`ere capacity improving a previous 
result in \cite{BJZ}.

\begin{prop} Let $\Om \Sub \C^n$ be a hyperconvex domain. Then for any Borel subset $E \sub
 \Om$, we have that
\begin{equation}
V_{2 n} (E) \leq c_n \delta _\Om^{2 n} \left(1 + Cap_{\Om} (E)^{- 1 \slash n} \right)^{n - 1} \exp \left(-  2 n \ Cap_{\Om} (E)^{- 1 \slash n}\right),
\label{eq:V-Cap}
\end{equation}
where  $\delta _\Om := diam (\Om)$ is the diameter of $\Om$ and $c_n$ is the constant defined by (\ref{eq:c_n}).
\end{prop}
\demo We first  assume that $E\Sub \Om $. Then its  \psh relative extremal function $h_E^* \in \mc E_0 (\Om)$. Therefore applying the last corollary, we obtain
$$
V_{2 n} (E) \leq  V_{2 n}  \left(\{ h_E^*  \leq - 1\}\right)   \leq  c_n \delta _\Om^{2 n} \left(1 +  \mu^{- 1} \right)^{n - 1} \exp \left(-  2 n \mu^{- 1}\right), 
$$
where $\mu^n = \int_{\Om}(dd^c h_E^*)^n$. Then the estimate of the theorem follows since 
 $\int_{\Om}(dd^c h_E^*)^n = Cap (E,\Om)$ by \cite{BT}.

Now assume that $Cap_{\Om} (E) < + \infty$. Then approximating $E$ by a non decreasing sequence of relatively compact subsets of $\Om$, it follows from continuity properties of the \MA operator in $\mc F (\Om)$ that
 $h_{E}^* \in \mc F (\Om)$ and the formula 
$\int_{\Om}(dd^c h_E^*)^n = Cap_{\Om} (E)$ still holds in this case. The proof of the inequality follows then in the same way.
\fin

Observe that actually the estimates (\ref{eq:V-Est}) and (\ref{eq:V-Cap}) are equivalent since for a function $\f \in \mc F (\Om)$ we know that
$Cap (\{\f \leq - s\} \leq s^{- n} \int_{\Om} (dd^c \f)^n$ (see \cite{CKZ}).

\medskip

\noindent{\bf Remark :}  Let $\Om \Sub \C^n$ be a hyperconvex domain such that $\Om \cap \R^n \neq \es$. Then the same method can be used to prove an estimate of the $n-$dimensional volume of Borel subsets of $ \Om \cap \R^n$ in terms of their capacity with respect to $\Om$. Namely, if $K \sub \Om \cap \R^n$ is a Borel subset, then its $n-$dimensional volume satisfies the inequality
$$
V_{n} (K) \leq b_n \delta _\Om^{2 n}  \left(1 + Cap_{\Om} (K)^{- 1 \slash n} \right)^{n - 1} \exp \left(-   n \ Cap_{\Om} (K)^{- 1 \slash n}\right),
$$
where $b_n > 0$ is a uniform constant which can be made explicit.
The proof uses induction as before and the following real version of  the inequality (\ref{eq:P}) (see \cite{Ra}): if $K \sub [- 1, + 1]$ is a real compact set of lenght  $V_1 (K)$ and logarithmic capacity $c (K)$ then
$$
V_1 (K) \slash 4  \leq c (K)  \leq 2 \exp (- 1 \slash Cap_{\D} (K)).
$$
(See \cite{BJZ} where such kind of estimates were considered).

\subsection{Integral estimates for capacity of slices}

 Given  a Borel subset $E \sub \C^{n + m} = \C^n \times \C^m,$ we define its 
 $n-$dimensional slices as follows. For a given $\zeta \in \C^m$ we  
 set
 $$ E_{\zeta} := \{z \in \C^{n} ; (z,\zeta) \in E \}.$$ 
 It is easy to see that if $E$ is pluripolar then its slices $E_{\zeta}$  
 are pluripolar sets in $\C^{n}$ except for a pluripolar set of $\zeta$'s in $\C^m$. The converse is not true as the following example of Kiselman \cite{Ki1}
$$ S := \{(z,w) \in \C^2 ;  Im (z+w^2) =Re (z+w+w^2) = 0\},$$
shows. Indeed $S$ is a smooth totally real analytic $2-$manifold in $\C^2$ whose intersection with any complex line consists of at most $4$ points.

 Here we want to give a quantitative estimate  in terms of \MA capacity of the size of  the slices of a Borel set.

 \begin{prop}  Let $\Om \sub \C^n$ and $D \Sub \C^m$ be two hyperconvex domains and $\Tilde \Om := \Om \times D.$  Assume that $E   
 \sub \tilde \Om$  is a Borel subset such that $Cap_{\tilde \Om} (E) < + \infty$. Then for any real number $p > 0$ and any we have that
 \begin{equation}
 \int_D \left(\mrm{Cap}_{\Om} (E_{\zeta})\right)^p d V_{2 m}(\zeta) \leq  
 (4 \pi)^m \  \de_D^{2 m} \ p^m \ 2^{- m p } \Gamma (p)^m Cap_{\tilde \Om} (E)^{p},
 \label{eq:1} 
 \end{equation}
 where $\Gamma$ is the Euler function defined by $\Gamma (p) :=  
 \int_0^{+ \infty} t^{p - 1} e^{- t} d t$. 
 \end{prop}
 \demo As in the proof of Proposition 6.1, if
$h_E^*$ is the \psh extremal function of the condenser $(E,\tilde \Om)$, we see that $h_E^* \in \mc F (\tilde \Om)$ if and only if $Cap_{\tilde \Om} (E) < + \infty.$ 

First  suppose  that $m = 1$ i.e. $D = \D \sub \C$ is the unit disc in $\C$. For each $\zeta \in \D,$ let $h_{E_{\zeta}}^*$ be the \psh extremal function of the condenser $(E_{\zeta},\Om)$. It follows from the definitions that for any $\zeta \in \D$,
 the partial function $h_E^* (\cdot,\zeta)$ satisfies the inequality
 $$h_E^* (\cdot,\zeta) \leq h_{E_{\zeta}}^*, \ \text{on} \ \Om.$$
 Moreover by Theorem 3.1, for almost all $\zeta \in \D$, these functions are in $\mc E_1 (\Om)$ and then by the proof of Lemma 2.1, we have
 \begin{eqnarray*}
 Cap_{\Om} (E_{\zeta})
  & = & \int_{\Om} (- h_{E_{\zeta}}^*) (dd^c h_{E_{\zeta}}^*)^{n} \\
  & \leq & \int_{\Om} (- h_E^* (\cdot,\zeta)) dd^c(h_E^* (\cdot,\zeta)^{n} \\
  & =:& - u (\zeta),
  \end{eqnarray*}
  where $u (\zeta)$ is  precisely the partial energy of the function  $h_{E}^* \in \mc F (\tilde \Om)$.
  By Theorem 3.1, we also have that $u \in \mc F (D)$ and 
    $$\int_{\D} dd^c u = \int_{\tilde \Om} (dd^c h_E^*)^{n + 1} = Cap_{\tilde \Om} (E).$$
  Now assume that $\chi: \R^+  \lra \R^+$ is an increasing function such that $\chi (0) = 0$, then applying (\ref{eq:VMA}) in the one dimensional case, we obtain
 \begin{eqnarray*}
  \int_{\D} \chi \left(Cap_{\Om} (E_{\zeta})\right) d V_2 (\zeta)
 & = &\int_0^{+ \infty} \chi' (s) V_2 (\{u \leq - s\} d s \\
 & \leq & 4 \pi  \int_0^{+ \infty} \chi' (s) \exp \left(- 2 s \slash Cap_{\tilde \Om} (E)\right) d s.
  \end{eqnarray*}
 If $\chi (t) = t^p$ then it follows that
 $$
 \int_{\D} \left(Cap_{\Om} (E_{\zeta})\right)^p d V_2 (\zeta)
 \leq 4 \pi p \int_0^{+ \infty} s^{p - 1} \exp \left(- 2 s \slash Cap_{\tilde \Om} 
 (E)\right) d s.
 $$
 Setting $t = 2 s \slash Cap_{\tilde \Om} (E)$, we obtain
 $$\int_{\D}  \left(Cap_{\Om} (E_{\zeta})\right)^p d V_2 (\zeta)
 \leq 4 \pi  p 2^{-p}  Cap_{\tilde \Om} (E)^{p} \int_0^{+ \infty} t^{p - 1} e^{- t}d t.
 $$
  Now if $D = \D^2 \Sub \C^2$ is is the unit polydisc , we can iterate the previous inequality. Observe that for any $(\zeta = (\zeta_1,\zeta_2) \in \D^2,$ we have that  $E_{\zeta_1} \sub \Om \times \D$ and
$E_{\zeta} = (E_{\zeta_1})_{\zeta_2}.$

 Therefore using the  previous estimate twice, we get
\begin{eqnarray*}
 \int_{\D^2} \left(Cap_{\Om} (E_{\zeta})\right)^p d V_4(\zeta_1,\zeta_2)
&=& \int_{\D} d V_2 (\zeta_1)  \int_{\D}  \left(Cap_{\Om} ((E_{\zeta_1})_{\zeta_2})\right)^p d V_2(\zeta_2) \\
&\leq & 4 \pi p \ 2^{- p } \Gamma (p) \int_{\D}  Cap_{ \Om \times \D} (E_{\zeta_1})^{p} d V_2 (\zeta_1) \\
& \leq & (4 \pi)^2 p^2 \ 2^{- 2 p } \Gamma (p)^2 Cap_{ \Om \times \D^2} (E)^{p}, \\
\end{eqnarray*} 

Now for $m \geq 3$ we obtain by induction on $m$ .

$$ \int_{\D^m} \left(Cap_{\Om} (E_{\zeta})\right)^p d V_{2 m} (\zeta)  
\leq  (4 \pi)^m p^m \ 2^{- 2 m p } \Gamma (p)^m Cap_{ \Om \times \D^m} (E)^{p}.$$
In the general case we can always assume that  $D \sub \D^n$ and  then $Cap_{\Om \times \D^n} (E) \leq Cap_{\Om \times D} (E)$ and the required estimate follows. 
  \fin
  
\subsection{A  local  transcendental inequality}

Here we want to give a  transcendental version of the local algebra inequality of  Corti \cite{Co} in dimension $2$, de Fernex, Ein and Musta\c{t}\v{a} \cite{DEM2} in higher dimensions.

Let us first recall the definition of complex integrability exponents introduced by  Demailly and  Koll\'{a}r \cite{DK}.
Let $\f$ be a plurisubharmonic function on an open set $\Om \sub \C^n$ and  $a \in \Om$.
 We define the complex singularity exponent of $\f$ at the point $a$ to be
$$c_a (\f) := \sup \{c > 0 ; \exists \ \ U \  \text{neighbourhood of } a, \ \exp (- 2 c \f) \in L^1_{loc}(U)\}.$$ 
 By Skoda's theorem \cite{Sk} we know that
$$\frac{1}{\nu_a (\f)} \leq c_a (\f) \leq \frac{n}{\nu_a (\f)},$$
where $\nu_a (\f)$ is the Lelong number of $\f$ at the point $a$ defined by the formula
$$
\nu_a (\f) := \sup \{\nu > 0 \ ;   \f (z) \leq \nu \log \vert z - a\vert + O (1), \ 0 <  \vert z - a\vert << 1\}.
$$ 
Our Theorem 5.4 can be rephrased in the following way.
\begin{prop}
Let $\f \in \mc E (\Om)$, then for any $a \in \Om$, we have 
$$ c_a (\f) \geq \frac{n}{\mu_a (\f)},$$
where $\mu_a (\f)$ is defined by the formula
$$ \mu_a (\f)^n := \int_{\{a\}} (dd^ c \f)^n.$$
Therefore
$$ \frac{n}{\mu_a (\f)} \leq  c_a (\f) \leq \frac{n}{\nu_a (\f)}.$$
\end{prop}
 As pointed out by Demailly \cite{De3}, this inequality implies an important inequality between two algebraic invariants associated to an ideal $\mc I$ of germs of holomorphic functions  with an isolated singularity at the origin in $\C^n$. Namely, if the ideal $\mc I$ is generated by the holomorphic germs $g_1, \cdots, g_N$ near the origin then its log canonical threshold at the origin is defined to be $lc (\mc I) := c_0 (\f),$ where $\f :=  \log \left(\sum_j \vert g_j\vert^2\right)$ (\cite{DK}).  There is another numerical  invariant $e (\mc I)$, called the Hilbert-Samuel multiplicity of the ideal $\mc I$ (see \cite{DEM1} for the definition), which turns out to be equal to 
$\mu_0 (\f)^n$ (\cite{De3}).

Thus our last result combined with Lemma 2.1 in \cite{De3} implies the following result from local algebra due to Corti \cite{Co} in dimension 2 and de Ferneque, Ein and
 Musta\c{t}\v{a} \cite{DEM1} in higher dimensions.
\begin{cor} Let $\mc I$ be an ideal as above. Then we have that
$$
 lc (\mc I)  \geq \frac{n}{\left(e (\mc I)\right)^{1 \slash n}}.
$$
\end{cor}

\newpage
\bigskip
\noindent Per {\AA}hag \\
Department of Mathematics and Mathematical Statistics \\
Ume\aa \ University \\
SE-901 87 Ume\aa \\
Sweden
\vskip 0.3 cm
\noindent Urban Cegrell \\
Department of Mathematics and Mathematical Statistics \\
Ume\aa \ University \\
SE-901 87 Ume\aa \\
Sweden

\vskip 0.3 cm
\noindent S\l awomir Ko\l odziej \\
Jagiellonian University \\
Institute of Mathematics \\
Reymonta 4, 30-059 Krak\'ow, Poland
\vskip 0.3 cm

\vskip 0.3 cm
\noindent Ph\d{a}m Hoang Hiep  \\
Department of Mathematics \\
Truong Dai Hoc Su Pham, \\
136 Xuan Thuy, Cau Giay, Hanoi, Vietnam
\vskip 0.3 cm

\noindent Ahmed Zeriahi \\
Universit\'e Paul Sabatier \\
Institut de Math\'ematiques \\
118 Route de Narbonne \\
31062 Toulouse cedex, France

\end{document}